\documentclass[12pt] {article}
\usepackage{amsmath,amsthm}
\usepackage{amssymb,latexsym}
\usepackage{graphicx}
\usepackage{amssymb}
\usepackage{amsfonts}
\usepackage{amscd}

\usepackage{amscd}
\title{Quaternionic plurisubharmonic functions and their applications to convexity.}
\date{}
\author{ Semyon Alesker \footnote{Partially supported by ISF grant 1369/04.}
\\  { \normalsize Department of Mathematics, Tel Aviv University, Ramat Aviv}
 \\  { \normalsize 69978 Tel Aviv,
Israel }
\\ {\normalsize e-mail: semyon@post.tau.ac.il}}

\def\eps{\varepsilon}

\def\ome{\omega}
\def\Ome{\Omega}
\def\lam{\lambda}
\def\Lam{\Lambda}

\def\to{\rightarrow}

\def\pt{\partial}
\def\d{\det}
\def\db{\frac{\partial}{\partial \bar q}}
\def\dq{\frac{\partial}{\partial  q}}

\def\RR{\mathbb{R}}
\def\CC{\mathbb{C}}

\def\HH{\mathbb{H}}

\def\duq{\frac{\partial ^2 u}{\partial\bar q_i \partial  q_j}}

\swapnumbers
\newtheorem{theorem}{Theorem}[section]
\newtheorem{corollary}[theorem]{Corollary}

\newtheorem{lemma}[theorem]{Lemma}
\newtheorem{proposition}[theorem]{Proposition}
\newtheorem{claim}[theorem]{Claim}
\theoremstyle{definition}
\newtheorem{example}[theorem]{Example}
\newtheorem{definition}[theorem]{Definition}
\newtheorem{remark}[theorem]{Remark}
\theoremstyle{conjecture}

\theoremstyle{proposition-definition}
\newtheorem{proposition-definition}[theorem]{Proposition-Definition}

 \def\ck{{\cal K}}

\begin{document}
\maketitle
\begin{quote}
\emph{ Dedicated to Professor Victor Abramovich Zalgaller in
occasion of his 85-th birthday.}
\end{quote}
\setcounter{section}{-1}
\begin{abstract}
The goal of this article is to present a survey of the recent
theory of plurisubharmonic functions of quaternionic variables,
and its applications to theory of valuations on convex sets and
HKT-geometry (HyperK\"ahler with Torsion). The exposition follows
the articles \cite{alesker1}, \cite{alesker2},
\cite{alesker-adv-05} by the author and \cite{alesker-verbitsky}
by M. Verbitsky and the author.
\end{abstract}
\tableofcontents
\section{Introduction.}\label{introduction}
The goal of this article is to present a survey of the recent
theory of plurisubharmonic functions of quaternionic variables,
and its applications to theory of valuations on convex sets and
HKT-geometry (HyperK\"ahler with Torsion). The exposition follows
the articles \cite{alesker1}, \cite{alesker2},
\cite{alesker-adv-05} by the author and \cite{alesker-verbitsky}
by M. Verbitsky and the author.

We will denote by $\HH$ the (non-commutative) field of
quaternions. The notion of quaternionic plurisubharmonic function
on the flat space $\HH^n$ was introduced by the author in
\cite{alesker1} and independently by G. Henkin \cite{henkin}. This
notion is a quaternionic analogue of convex functions on $\RR^n$
and complex plurisubharmonic functions on $\CC^n$, see Definition
\ref{def-qpsh} below. On one hand, this class of functions obeys
many analytical properties analogous to those of convex and
complex plurisubharmonic functions. On the other hand, these
properties reflect rather different geometric structures behind
them. This will be illustrated below on applications to convexity
and HKT-geometry.

Let us start with some analytical properties of quaternionic
plurisubharmonic functions. The author has proved in
\cite{alesker1} a quaternionic analogue of the A.D. Aleksandrov
\cite{aleksandrov2} and Chern-Levine-Nirenberg
\cite{chern-levine-nirenberg} theorems (see Theorems \ref{cln1},
\ref{cln2} in Section \ref{psh} below). It is worthwhile to remind
these classical results now. The Aleksandrov theorem says that if
a sequence $\{f_N\}$ of {\itshape convex} functions converges
uniformly on compact subsets to a function $f$, then $f$ is convex
and
$$\det\left(\frac{\pt^2f_N}{\pt x_i\pt x_j}\right)\overset{w}{\to}\det\left(\frac{\pt^2f}{\pt x_i\pt
x_j}\right)$$ weakly in the sense of measures (note that the
expression $\det\left(\frac{\pt^2u}{\pt x_i\pt x_j}\right)$ for a
convex function $u$ is understood in a generalized sense as
explained in the quaternionic situation in Theorem \ref{cln1} of
Section \ref{psh} below).

The Chern-Levine-Nirenberg theorem (in fact, in a slightly weaker
form) says that if a sequence $\{f_N\}$ of {\itshape continuous
complex plurisubharmonic} functions converges uniformly on compact
subsets to a function $f$, then $f$ is continuous complex
plurisubharmonic and
$$\det\left(\frac{\pt^2f_N}{\pt z_i\pt \bar
z_j}\right)\overset{w}{\to}\det\left(\frac{\pt^2f}{\pt z_i\pt \bar
z_j}\right)$$ weakly in the sense of measures (again the
expressions $\det\left(\frac{\pt^2u}{\pt z_i\pt \bar z_j}\right)$
is understood in a generalized sense).

The statement of the quaternionic analogue of the above theorems
requires analogues of complex operators $\frac{\pt}{\pt \bar z},\,
\frac{\pt}{\pt z}$ and the notion of determinant of quaternionic
matrices. The former notion is called sometimes Dirac operators
$\frac{\pt}{\pt\bar q},\, \frac{\pt}{\pt q}$; it is discussed in
Section \ref{dirac}. The latter notion of quaternionic
determinants is discussed in Section \ref{linear_algebra} where we
discuss in detail the Moore determinant of hyperhermitian ( =
quaternionic hermitian) matrices. This quaternionic result is used
in applications to theory of valuations on convex sets (Theorem
\ref{q-kazgen} in Section \ref{valuations}).

Another important for applications in valuation theory result is
Theorem \ref{blocki} (proved in \cite{alesker-adv-05}) which is a
quaternionic version of B\l ocki's theorem \cite{blocki} for
complex plurisubharmonic functions.

In Section \ref{valuations} we discuss in more details
applications of the above results to this theory of valuations on
convex sets. Let us remind basic notions of this theory referring
for further information to the surveys by McMullen
\cite{mcmullen-survey} and McMullen and Schneider
\cite{mcmullen-schneider}. Let $V$ be a finite dimensional real
vector space. Let ${\cal K}(V)$ denote the class of all non-empty
convex compact subsets of $V$.

\begin{definition}
(1) A function $\phi :{\cal K}(V) \to \CC$ is called a valuation
if for any $K_1, \, K_2 \in {\cal K}(V)$ such that their union is
also convex one has
$$\phi(K_1 \cup K_2)= \phi(K_1) +\phi(K_2) -\phi(K_1 \cap K_2).$$

(2) A valuation $\phi$ is called continuous if it is continuous
with respect the Hausdorff metric on ${\cal K}(V)$.
\end{definition}

Recall that the Hausdorff metric $d_H$ on ${\cal K}(V)$ depends on
a choice of a Euclidean metric on $V$ and it is defined as
follows: $d_H(A,B):=\inf\{ \eps >0|A\subset (B)_\eps \mbox{ and }
B\subset (A)_\eps\}$ where $(U)_\eps$ denotes the
$\eps$-neighborhood of a set $U$. Then ${\cal K} (V)$ becomes a
locally compact space, and the topology on $\ck(V)$ induced by the
Hausdorff metric does not depend on a choice of the Euclidean
metric on $V$.

The theory of valuations has numerous applications in convexity
and integral geometry (see e.g. \cite{alesker-jdg},
\cite{hadwiger-book}, \cite{klain-rota}, \cite{schneider-book}).
Let us remind some basic examples of translation invariant
continuous valuations.
\begin{example}\label{ex-valuations}
(1) A Lebesgue measure $vol$ on $V$ is a translation invariant
continuous valuation.

(2) The Euler characteristic $\chi$ is a translation invariant
continuous valuation. (Recall that $\chi (K)=1$ for any $K\in
\ck(V)$.)

(3) Denote $m:=\dim V$. Fix $k=1,\dots,m$. Fix
$A_1,\dots,A_{m-k}\in \ck(V)$. Then the mixed volume $$K\mapsto
V(K[k],A_1,\dots,A_{m-k})$$ is a translation invariant continuous
valuations. (For the notion of mixed volume and its properties see
e.g. the books by Burago-Zalgaller \cite{burago-zalgaller} and
Schneider \cite{schneider-book}.)
\end{example}
It was conjectured by P. McMullen \cite{mcmullen-80} and proved by
the author \cite{alesker-gafa-01} that the linear combinations of
the mixed volumes as in Example \ref{ex-valuations} (3) above are
dense in the space of all translation invariant continuous
valuations in the topology of uniform convergence on compact
subsets of $\ck(V)$.

Nevertheless there are other than mixed volumes non-trivial
constructions of translation invariant continuous valuations. In
this survey we will discuss two of them. They are closely related
to each other since they are based on the theory of complex and,
respectively, quaternionic plurisubharmonic functions. Their
relation to the mixed volume construction is not straightforward.
The possibility of approximating of these examples by the mixed
volumes follows from the solution of the McMullen's conjecture.

Let us agree on the notation. For $K\in \ck(V)$ one denotes by
$h_K\colon V^*\to \RR$ the supporting functional of $K$. Recall
that
$$h_K(y)=\sup\{y(x)|\, y\in K\}.$$

Let us describe the construction of valuations using complex
plurisubharmonic functions. Let us denote by $\Ome^{p,p}$ the
vector bundle over $\CC^n$ of $(p,p)$-forms. Let us denote by
$C_c(\CC^n,\Ome^{p,p})$ the space of continuous compactly
supported forms of type $(p,p)$ on $\CC^n$.
\begin{theorem}[\cite{alesker-adv-05}, Theorem 4.1.3]\label{kazgen-intro}
Fix $k=1,\dots,n$.  Fix $\psi\in C_0(\CC^n, \Omega^{n-k,n-k})$.
Then $K\mapsto \int_{\CC^n}(dd^ch_K)^k\wedge \psi$ defines a
continuous translation invariant valuation on $\ck (\CC^n)$.
\end{theorem}
This result in the above generality was proved by the author in
\cite{alesker-adv-05} using some known properties of complex
plurisubharmonic functions. This theorem contains two non-trivial
parts: the continuity and the valuation property of the above
functional. The former is a consequence of the
Chern-Levine-Nirenberg theorem \cite{chern-levine-nirenberg}, and
the latter is a consequence of the B\l ocki formula \cite{blocki}.
Note that the expressions of the form as in Theorem \ref{kazgen}
were considered first in the context of convexity probably by
Kazarnovski\u\i  \cite{kazarnovskii-81}, \cite{kazarnovskii-84}.
In \cite{alesker-adv-05} we have obtained quaternionic version of
the above construction; it is discussed in Section
\ref{valuations}, Theorem \ref{q-kazgen}.

In Section \ref{monge-ampere} we discuss two theorems on the
Dirichlet problem for quaternionic Monge-Amp\`ere equations
obtained by the author in \cite{alesker2}. They also have
classical real and complex analogues; we refer to Section
\ref{monge-ampere} for references.

In Section \ref{hypercomplex} we describe generalizations of some
of the definitions and results on quaternionic plurisubharmonic
functions to so called hypercomplex manifolds due to M. Verbitsky
and the author \cite{alesker-verbitsky}. This class of manifolds
contains, for instance, the flat spaces $\HH^n$ and hyperK\"ahler
manifolds. On such manifolds one defines quaternionic
plurisubharmonic functions and proves for them a version of the
Aleksandrov and Chern-Levine-Nirenberg theorems. Next it turns out
that $C^\infty$-smooth strictly plurisubharmonic functions on
hypercomplex manifolds admit a geometric interpretation as (local)
potentials of HKT-metrics; it was also obtained in
\cite{alesker-verbitsky}. Roughly put, an HKT-metrics on a
hypercomplex manifold is an $SU(2)$-invariant Riemannian metric
satisfying certain first order differential equations. These
metrics are analogous to K\"ahler metrics on complex manifolds.
The above mentioned interpretation of quaternionic
plurisubharmonic functions is analogous to the well known
interpretation of $C^\infty$-smooth complex strictly
plurisubharmonic functions on complex manifolds as (local)
potentials of K\"ahler metrics. Let us recall more explicitly the
last fact. Let $f$ be a $C^\infty$-smooth complex strictly
plurisubharmonic function on a complex manifold $M$.  Let us fix
on $M$ local complex coordinates. Then the matrix
$g:=\left(\frac{\pt^2f}{\pt z_i\pt \bar z_j}\right)$ defines a
K\"ahler metric on $M$. Vice versa, for any K\"ahler metric $g$ on
$M$, every point $z\in M$ has a neighborhood such that in this
neighborhood $g=\left(\frac{\pt^2f}{\pt z_i\pt \bar z_j}\right)$
for a $C^\infty$-smooth complex strictly plurisubharmonic function
$f$.

The organization of the article is clear from the table of
contents.

{\bf Acknowledgements.} I thank G. Henkin and M. Verbitsky for
very useful discussions.

\section{Quaternionic linear algebra.}\label{linear_algebra}
It is known that some of the standard results of linear algebra
over commutative fields can be generalized to  general
non-commutative fields, e.g. theory of dimension and basis of
vector spaces. However over non-commutative fields there is no
notion of determinant of matrices which would be as good as in the
commutative case. There is a notion of the Diedonn\'e determinant
(see e.g. \cite{artin}) which is good for some applications (e.g.
\cite{alesker1}). Gelfand and Retakh have developed a theory of
quasi-determinants over non-commutative fields which generalizes
in a sense probably all the known theories of non-commutative
determinants (in particular the Dieudonn\'e, super-, quantum-
determinants). We will not discuss this theory here, and refer to
the survey \cite{ggrw}.

Nevertheless over the field $\HH$ of quaternions there is a notion
of {\itshape Moore determinant} on the class of {\itshape
hyperhermitian} matrices discussed below. Hyperhermitian
quaternionic matrices are analogous to real symmetric and complex
hermitian matrices. The properties of the Moore determinant on
this class are very similar to the properties of the usual
determinant. It seems that any general identity or inequality
which is true for the determinant of real symmetric or complex
hermitian matrices should be true for the Moore determinants of
quaternionic hyperhermitian matrices. Among examples of such
results one can mention Sylvester criterion of positive
definiteness and Aleksandrov inequalities for mixed determinants
discussed below. For more information on quaternionic determinants
see the survey \cite{aslaksen} and references therein; for the
relation of the Moore determinant to the Gelfand-Retakh
quasideterminants see \cite{gelfand-retakh-wilson}.

In the rest of this section we discuss in more detail the notion
of hyperhermitian matrices and their Moore determinant.

\begin{definition}
  A {\itshape hyperhermitian semilinear form}
on $V$ is a map $ a:V \times V \to \HH$ satisfying the following
properties:

(a) $a$ is additive with respect to each argument;

(b) $a(x,y \cdot q)= a(x, y) \cdot q$ for any $x,y \in V$ and any
$q\in \HH$;

(c) $a(x,y)= \overline{a(y,x)}$ where $\bar q$ denotes the usual
conjugation of a quaternion $q$
\end{definition}
\begin{definition}
A quaternionic $(n\times n)$-matrix $A=(a_{ij})_{i,j=1}^n$ is
called {\itshape hyperhermitian} if $a_{ij}=\bar a_{ji}$.
\end{definition}
\begin{example}\label{ex-hyperherm} Let $V= \HH ^n$ be the standard coordinate space
considered as right vector space over $\HH$. Fix a {\itshape
hyperhermitian} $(n \times n)$-matrix $(a_{ij})_{i,j=1}^{n}$. For
$x=(x_1, \dots, x_n), \, y=(y_1, \dots, y_n)$ define
$$a(x,y) = \sum _{i,j} \bar x_i a_{ij} y_j$$
(note the order of the terms!). Then $a$ defines hyperhermitian
semilinear form on $V$.
\end{example}
In general one has the following standard claims.

\begin{claim} Fix a basis in a finite dimensional right quaternionic
vector space $V$. Then there is a natural bijection between
hyperhermitian semilinear forms on $V$ and $(n \times
n)$-hyperhermitian matrices.
\end{claim}
This bijection is in fact described in previous Example
\ref{ex-hyperherm}.

\begin{claim} Let $A$ be the matrix
of a given hyperhermitian form in a given basis. Let $C$ be
transition matrix from this basis to another one. Then the matrix
$A'$ of the given form in the new basis is equal to $C^* AC .$
\end{claim}

\begin{definition} A hyperhermitian semilinear form $a$
is called {\itshape positive definite} if $a(x,x)>0$ for any
non-zero vector $x$.
 \end{definition}

Let us fix on our quaternionic right vector space $V$ a positive
definite hyperhermitian form $( \cdot , \cdot )$. The space with
fixed such a form will be called {\itshape hyperhermitian} space.

For any quaternionic linear operator $\phi: V\to V$ in
hyperhermitian space $V$ one can define the adjoint operator $\phi
^* :V \to V$ in the usual way, i.e. $(\phi x,y)= (x, \phi ^* y)$
for any $x,y \in V$. Then if one fixes an orthonormal basis in the
space $V$ then the operator $\phi$ is selfadjoint if and only if
its matrix in this basis is hyperhermitian.

\begin{claim}
For any  selfadjoint operator in a hyperhermitian space there
exists an orthonormal basis such that its matrix in this basis is
diagonal and real.
\end{claim}
Now we are going to define the Moore determinant of
hyperhermitian matrices. The definition below is different from
the original one due to Moore \cite{moore} but equivalent to it.

First note that every hyperhermitian $(n \times n)$- matrix $A$
defines a hyperhermitian semilinear form on the coordinate space
$\HH ^n$ as explained in Example \ref{ex-hyperherm}. It also can
be considered as a {\itshape symmetric} bilinear form on $\RR
^{4n}$ (which is the realization of $\HH ^n$). Let us denote its
$(4n \times 4n)$- matrix by ${}^{\textbf{R}} A$. Let us consider
the entries of $A$ as formal variables (each quaternionic entry
corresponds to four commuting real variables). Then $det
({}^{\textbf{R}} A)$  is a homogeneous polynomial of degree $4n$
in $n(2n-1)$ real variables. Let us denote by $Id$ the identity
matrix. One has the following result which was rediscovered
several times (see \cite{aslaksen} for references).

\begin{theorem}
There exists a polynomial $P$ defined on the space of all
hyperhermitian $(n \times n)$-matrices such that for any
hyperhermitian $(n \times n)$-matrix $A$ one has
$det({}^{\textbf{R}} A)= P^4(A)$ and $P(Id)=1$. $P$ is defined
uniquely by these two properties. Furthermore $P$ is homogeneous
of degree $n$ and has integer coefficients.
\end{theorem}
Thus for any hyperhermitian matrix $A$ the value $P(A)$ is a real
number, and it is called the {\itshape Moore determinant} of the
matrix $A$. The explicit formula for the Moore determinant  was
given by Moore \cite{moore} (see also \cite{aslaksen} and Theorem
\ref{moore-formula} below). From now on the Moore determinant of a
matrix $A$ will be denoted by $det A$. This notation should not
cause any confusion with the usual determinant of real or complex
matrices due to part (i) of the next theorem which is probably a
folklore.

\begin{theorem}

(i) The Moore determinant of any complex hermitian matrix
considered as quaternionic hyperhermitian matrix is equal to its
usual determinant.

(ii) For any hyperhermitian matrix $A$ and any quaternionic matrix
$C$
$$det (C^*AC)= detA \cdot det(C^*C).$$
\end{theorem}
\begin{example}

(a) Let $A =diag(\lam_1, \dots, \lam _n)$ be a diagonal matrix
with real $\lam _i$'s. Then $A$ is hyperhermitian and the Moore
determinant $detA= \prod _i \lam_i$.

(b)  A general hyperhermitian $(2 \times 2)$-matrix $A$ has the
form
 $$ A=  \left[ \begin {array}{cc}
                     a&q\\
                \bar q&b\\
                \end{array} \right] ,$$
where $a,b \in \RR, \, q \in \HH$. Then $det A =ab - q \bar q (
=ab- \bar q q)$.
\end{example}

\begin{claim}
Let $A$ be  a non-negative (resp. positive) definite
hyperhermitian matrix. Then $det A \geq 0 \, (\mbox{ resp. } det A
>0)$.
\end{claim}

The following theorem is a quaternionic generalization of the
standard Sylvester criterion.
\begin{theorem}[Sylvester criterion, \cite{alesker1}]
A hyperhermitian $(n \times n)$- matrix $A$ is positive definite
if and only if the Moore determinants of all the left upper minors
of $A$ are positive.
\end{theorem}

Let us define now the mixed discriminant of hyperhermitian
matrices in analogy with the case of real symmetric matrices
studied by A.D. Aleksandrov \cite{aleksandrov1}.
\begin{definition}
Let $A_1, \dots ,A_n$ be hyperhermitian $(n \times n)$- matrices.
Consider the homogeneous polynomial in real variables $\lam _1
,\dots , \lam _n$ of degree $n$ equal to $det(\lam_1 A_1 + \dots +
\lam_n A_n)$. The coefficient of the monomial $\lam_1 \cdot \dots
\cdot \lam_n$ divided by $n!$ is called the {\itshape mixed
discriminant} of the matrices  $A_1, \dots ,A_n$, and it is
denoted by $\d(A_1, \dots ,A_n)$.
\end{definition}
Note that the mixed discriminant is symmetric with respect to all
variables, and linear with respect to each of them; recall that
the linearity with respect to say the first argument means that
$$\d (\lam A_1' +\mu A_1'', A_2, \dots, A_n )=
\lam \cdot \d( A_1', A_2, \dots, A_n ) + \mu \cdot \d(A_1'', A_2,
\dots, A_n )$$ for any {\itshape real} $\lam , \, \mu$. Note also
that $\d(A, \dots, A)=det A$. One has the following generalization
of Aleksandrov's inequalities for mixed discriminants
\cite{aleksandrov1}.
\begin{theorem}
(i) The mixed discriminant of positive (resp. non-negative)
definite matrices is positive (resp. non-negative).

(ii) Fix positive definite hyperhermitian $(n \times n)$- matrices
$A_1, \dots, A_{n-2}$. On the real linear space of hyperhermitian
$(n \times n)$- matrices consider the bilinear form $$B(X,Y):=
\d(X,Y, A_1, \dots, A_{n-2}).$$ Then $B$ is non-degenerate
quadratic form, and its signature has one plus and the rest are
minuses.
\end{theorem}
\begin{corollary}[Aleksandrov inequality, \cite{alesker1}]
 Let
$A_1, \dots, A_{n-1}$ be positive definite hyperhermitian $(n
\times n)$- matrices. Then for any hyperhermitian matrix $X$
$$
\d(A_1, \dots, A_{n-1}, X)^2 \geq \d(A_1, \dots, A_{n-1},A_{n-1})
\cdot \d(A_1, \dots , A_{n-2}, X,X) ,
$$
and the equality is satisfied if and only if the matrix $X$ is
proportional to $A_{n-1}$.
\end{corollary}

Finally let us give an explicit formula for the Moore determinant
(which was the original definition by Moore \cite{moore}). Let
$A=(a_{i,j})_{i,j=1}^n$ be a hyperhermitian $(n\times n)$-matrix.
Let $\sigma$ be a permutation of $\{1,\dots, n\}$. Write $\sigma$
as a product of disjoint cycles such that each cycle starts with
the smallest number. Since disjoint cycles commute we may write
$$\sigma=(k_{11}\dots k_{1 j_1})(k_{21}\dots k_{2 j_2})\dots
(k_{m1} \dots k_{m j_m})$$ where for each $i$ we have
$k_{i1}<k_{ij}$ for all $j>1$, and $k_{11}>k_{21}>\dots >k_{m1}$.
This expression is unique. Let $sgn(\sigma)$ is the parity of
$\sigma$. For the next result we refer to \cite{aslaksen} and
references therein.
\begin{theorem}\label{moore-formula}
The Moore determinant of $A$ is equal to
$$\det A= \sum _\sigma sgn (\sigma) a_{k_{11},k_{12}}\dots
a_{k_{1 j_1},k_{11}} a_{k_{21},k_{22}}\dots a_{k_{m j_m},
k_{m1}}$$ where the sum runs over all permutations $\sigma$.
\end{theorem}

\section{Dirac operators.}\label{dirac}
We will write a quaternion $q\in \HH$ in the standard form
$$q= t+ x\cdot i +y\cdot j+ z\cdot k ,$$
where $t,\, x,\, y,\, z$ are real numbers, and $i,\, j,\, k$
satisfy the usual relations
$$i^2=j^2=k^2=-1, \, ij=-ji=k,\, jk=-kj=i, \, ki=-ik=j.$$

The Dirac operator $\frac {\partial}{\partial \bar q}$ is defined
as follows. For any $\HH$-valued function $F$
$$\db F:=\frac{\partial F}{\partial  t}  +
i \frac{\partial F}{\partial  x} + j \frac{\partial F}{\partial y}
+ k \frac{\partial F}{\partial  z}.$$

Let us also define the operator $\dq$:
$$\dq F:=\overline{ \db \bar F}=
\frac{\partial F}{\partial  t}  -
 \frac{\partial F}{\partial x}  i-
 \frac{\partial F}{\partial  y} j-
\frac{\partial F}{\partial  z}  k.$$

In the case of several quaternionic variables, it is easy to see
that the operators $\frac{\pt}{\pt q_i}$ and $\frac{\pt}{\pt \bar
q_j}$ commute:
\begin{eqnarray}
\big[\frac{\pt}{\pt q_i},\frac{\pt}{\pt \bar q_j}\big]=0.
\end{eqnarray}
\begin{proposition}[\cite{alesker1}]
(i) Let $f:\HH ^n \to \HH$ be a smooth function. Then for any
$\HH$-linear transformation $A$ of $\HH ^n$ (as a right $\HH
$-vector space) one has the identities
$$ \left( \frac {\partial ^2 f(Aq)}{\partial \bar q_i \partial q_j} \right)
=A^* \left(\frac {\partial ^2 f}{\partial \bar q_i \partial
q_j}(Aq) \right)A .$$

(ii) If, in addition, $f$ is real valued then for any $\HH$-linear
transformation $A$ of $\HH ^n$ and any quaternion $a$ with $|a|=1$
$$ \left( \frac {\partial ^2 f(A(q \cdot a))}{\partial \bar q_i \partial q_j} \right)
=A^* \left(\frac {\partial ^2 f}{\partial \bar q_i \partial
q_j}(A(q\cdot a)) \right) A .$$
\end{proposition}

\section{Plurisubharmonic functions of quaternionic
variables.}\label{psh}First we introduce the class of quaternionic
plurisubharmonic functions on $\HH^n$ following \cite{alesker1}.
Note that this notion was also introduced independently by G.
Henkin \cite{henkin}.
 Let $\Omega $ be an open subset of $\HH ^n$.
\begin{definition}\label{def-qpsh}
A real valued function $u: \Omega \to \RR$ is called quaternionic
plurisubharmonic if it is upper semi-continuous and its
restriction to any right {\itshape quaternionic} line is
subharmonic.
\end{definition}
Recall that upper semi-continuity means that $u(x_0)\geq
\underset{x\to x_0}{\limsup u(x)}$ for any $x_0\in \Omega$. We
will denote by $P(\Omega)$ the class of plurisubharmonic functions
in the open subset $\Omega$.

Also we will call a $C^2$-smooth function $u:\Omega\to \RR$ to be
{\itshape strictly plurisubharmonic} if its restriction to any
right quaternionic line is strictly subharmonic (i.e. the
Laplacian is strictly positive).

Before we state the next proposition let us observe that if a
smooth function $f$ is real valued then the matrix
$(\frac{\partial ^2 f}{\partial\bar q_i \partial  q_j})(q)$ is
hyperhermitian.
\begin{proposition}[\cite{alesker1}, Prop. 2.1.6]\label{psh-char}
A real valued twice continuously differentiable function $f$ on an
open subset $\Omega \subset \HH ^n$ is quaternionic
plurisubharmonic (reps. strictly plurisubharmonic) if and only if
at every point $q \in \Omega$ the matrix $(\frac{\partial ^2
f}{\partial\bar q_i
\partial q_j})(q)$ is non-negative definite (resp. positive definite).
\end{proposition}
\begin{remark}
Proposition \ref{psh-char} is completely analogous to
characterization of smooth convex function as functions with
non-negative definite Hessian $\big(\frac{\pt^2 f}{\pt x_i\pt
x_j}\big)$, and smooth complex plurisubharmonic functions as
functions with non-negative definite complex Hessian
$\big(\frac{\pt^2 f}{\pt z_i\pt \bar z_j}\big)$.
\end{remark}

\begin{theorem}[\cite{alesker1}]\label{cln1}
For any function $u\in C(\Omega) \cap P(\Omega)$ one can uniquely
define a non-negative measure   $\det (\duq)$ which is uniquely
characterized by the following two properties:
\newline
(1) if $u\in C^2(\Omega)$ then it has the obvious meaning;
\newline
(2) if $u_N\to u$ uniformly on compact subsets in $\Omega$, and
$u_N,\, u\in C(\Omega)\cap P(\Omega)$, then
$$\det (\frac{\partial ^2 u_N}{\partial\bar q_i \partial q_j})\overset
{w}{\to} \det(\duq),$$ where the convergence of measures in
understood in the weak sense.
\end{theorem}

\begin{remark} (1) It is easy to see that if $u_N\to u$ uniformly on
compact subsets, and $u_N\in C(\Omega)\cap P(\Omega)$ then $u\in
C(\Omega)\cap P(\Omega)$.

(2) Note that the real analogue of this result was proved by A.D.
Aleksandrov \cite{aleksandrov2}, and the complex analogue by
Chern, Levine, and Nirenberg \cite{chern-levine-nirenberg}.
\end{remark}
We will need a refinement of Theorem \ref{cln1} which was proved
by the author in \cite{alesker-adv-05} in somewhat different
notation. We denote for brevity $\pt^2u:=\left(\duq\right)$.

\begin{theorem}[\cite{alesker-adv-05}]\label{cln2}
Let $\Ome\subset \HH^n$ be an open subset. Fix $k=1,\dots,n$. Let
$\{u^{(i)}_N\}_{N=1}^\infty$, $1\leq i\leq k$, be sequences in
$P(\Ome)\cap C(\Ome)$. Let $V^{(1)},\dots,V^{(n-k)}$ be continuous
functions on $\Ome$ with values in the space of $(n\times
n)$-hyperhermitian matrices. Assume that for every $i=1,\dots,k$
$$u^{(i)}_N\to u^{(i)} \mbox{ as } N\to \infty$$ uniformly on compact subsets. Then
$u^{(i)}\in P(\Ome)\cap C(\Ome)$, and
$$\det\left(\pt^2 u^{(1)}_N,\dots,\pt^2
u^{(k)}_N,V^{(1)},\dots,V^{(n-k)}\right)\overset{w}{\to}
\det\left(\pt^2 u^{(1)},\dots,\pt^2
u^{(k)},V^{(1)},\dots,V^{(n-k)}\right)$$ weakly in the sense of
measures (where we used the notion of mixed determinant).
\end{theorem}
For hyperhermitian matrices $A,B_1,\dots,B_{n-k}$ let us denote
$$\det(A[k],B_1,\dots,B_{n-k}):= \det(\underset{k \mbox{
times}}{\underbrace{A,\dots,A}},B_1,\dots,B_{n-k}).$$

Note that the maximum of two plurisubharmonic functions is again
plurisubharmonic.
\begin{theorem}[\cite{alesker-adv-05}]\label{blocki}
Let $\Ome\subset \HH^n$ be an open subset. Fix $k=1,\dots,n$. Let
$f,g\in P(\Ome)\cap C(\Ome)$. Assume that $\min\{f,g\}\in
P(\Ome)\cap C(\Ome)$. Let $V^{(1)},\dots,V^{(n-k)}$ be continuous
functions on $\Ome$ with values in the space of $(n\times
n)$-hyperhermitian matrices. Then
\begin{eqnarray*}
\det(\pt^2(\max\{f,g\})[k],V^{(1)},\dots,V^{(n-k)})=\\
\det(\pt^2f[k],V^{(1)},\dots,V^{(n-k)})+\det(\pt^2g[k],
V^{(1)},\dots,V^{(n-k)})-\\
\det(\pt^2(\min\{f,g\})[k],V^{(1)},\dots,V^{(n-k)}).
\end{eqnarray*}
\end{theorem}
\begin{remark}
Theorem \ref{blocki} was proved in \cite{alesker-adv-05} as a
consequence of more precise result, Theorem 3.2.1 there, which is
a quaternionic version of a result by B\l ocki \cite{blocki} for
complex plurisubharmonic functions.
\end{remark}

\section{Applications to valuation theory.}\label{valuations}
Let us discuss the applications of the above results to the theory
of valuations on convex sets. The necessary definitions of this
theory were reminded in Introduction. Here we would like to recall
Theorem \ref{kazgen-intro}.
\begin{theorem}[\cite{alesker-adv-05}, Theorem 4.1.3]\label{kazgen}
Fix $k=1,\dots,n$.  Fix $\psi\in C_0(\CC^n, \Omega^{n-k,n-k})$.
Then $K\mapsto \int_{\CC^n}(dd^ch_K)^k\wedge \psi$ defines a
continuous translation invariant valuation on $\ck (\CC^n)$.
\end{theorem}
This result in the above generality was proved by the author in
\cite{alesker-adv-05} using some known properties of complex
plurisubharmonic functions. This theorem contains two non-trivial
parts: the continuity and the valuation property of the above
functional. The continuity is a consequence of the
Chern-Levine-Nirenberg theorem \cite{chern-levine-nirenberg} and
the fact that a sequence $\{K_N\}$ of convex compact sets
converges in the Hausdorff metric to a convex compact set $K$ if
and only if $h_{K_N}\to h_K$ uniformly on compact subsets. The
valuation property is a consequence of the B\l ocki formula
\cite{blocki} combined with the facts that
\begin{eqnarray*}
h_{K_1\cup K_2}=\max\{h_{K_1},h_{K_2}\} \mbox{ when } K_1\cup K_2
\mbox{ is convex, and}\\
h_{K_1\cap K_2}=\min\{h_{K_1},h_{K_2}\}.
\end{eqnarray*}

Let us discuss the quaternionic version of Theorem \ref{kazgen}.
In order to simplify the exposition we will state the result in a
simple minded form. A better way was discussed in
\cite{alesker-adv-05} where quaternionic analogues of
$(p,p)$-forms were introduced.
\begin{theorem}[\cite{alesker-adv-05}, Theorem
4.2.1]\label{q-kazgen} Fix $k=1,\dots,n$. Let $\psi_0$ be a
continuous compactly supported real valued function on $\HH^{n*}$.
Let $V^{(1)},\dots,V^{(n-k)}$ be continuous compactly supported
functions on $\HH^{n*}$ with values in the space of $(n\times
n)$-hyperhermitian matrices. Then the functional
$$K\mapsto \int_{\HH^{n*}}
\det(\pt^2h_K[k],V^{(1)},\dots,V^{(n-k)})\cdot \psi_0\cdot dvol$$
is a translation invariant continuous valuation.
\end{theorem}
Similarly to the complex case, the continuity is a consequence of
Theorem \ref{cln2}, and the valuation property is a consequence of
Theorem \ref{blocki}.
\section{Quaternionic Monge-Amp\`ere
equations.}\label{monge-ampere} In this section we will discuss
some results on quaternionic Monge-Amp\`ere equations following
\cite{alesker2}. They have real and complex analogues; the
references will be given below.
\begin{definition}
  An open bounded domain $\Omega\subset
\HH^n$ with a smooth boundary $\partial \Omega$ is called strictly
pseudoconvex if for every point $z_0\in \partial \Omega$ there
exists a neighborhood ${\cal O}$ and a smooth quaternionic
strictly plurisubharmonic function $h$ on ${\cal O}$ such that
$\Omega \cap {\cal O}= \{ h<0 \}$, $ h(z_0)=0$, and $\nabla h(z_0)
\ne 0$.
\end{definition}
Let $B$ denote the unit Euclidean ball in $\HH^n$.
\begin{theorem}[\cite{alesker2}, Theorem 0.1.4]\label{mareg}
Let $f\in C^{\infty}(\bar B), \, f>0$. Let $\phi \in
C^{\infty}(\partial B)$. There exists unique function $u\in
C^{\infty}(\bar B)$ which is quaternionic plurisubharmonic in $B$
and which is a solution of the Dirichlet problem
$$\det (\frac{\partial ^2 u}{\partial\bar q_i \partial q_j})=f \mbox{ in } B,$$
$$u|_{\partial B}= \phi.$$
\end{theorem}
\begin{remark}
(1) It is natural to expect that Theorem \ref{mareg} is true for a
larger class of domains, say for all bounded strictly pseudoconvex
domains.

(2)The real version of Theorem \ref{mareg} was proved for
arbitrary bounded strictly convex domains in $\RR^n$ by
Caffarelli, Nirenberg, and Spruck \cite{CNS}. The complex version
of it was proved for arbitrary bounded strictly pseudoconvex
domains in $\CC^n$ by Caffarelli, Kohn, Nirenberg, and Spruck
\cite{CKNS} and Krylov \cite{krylov2}. Our method is a
modification of the method of the
 paper \cite{CKNS}. Note also that in the case $n=1$ the
problem is reduced to the classical Dirichlet problem for the
Laplacian in $\RR^4$ (which is a linear problem); it was solved in
XIX century. Note also that {\itshape interior} regularity of the
solution of the Dirichlet problem for real Monge-Amp\`ere
equations was proved earlier by A. Pogorelov, and the proof was
briefly described  in \cite{pogorelov1}-\cite{pogorelov3}. The
complete proof was published in \cite{pogorelov} and
\cite{cheng-yau1}, \cite{cheng-yau2}.\end{remark}
\begin{theorem}[\cite{alesker2}, Theorem 0.1.3]
Let $\Omega \subset \HH ^n$ be a bounded quaternionic strictly
pseudoconvex domain. Let $f\in C(\bar \Omega), \, f \geq 0$. Let
$\phi \in C(\partial \Omega)$. Then there exists unique function
$u\in C(\bar \Omega)$ which is plurisubharmonic in $\Omega$ and
such that
$$\det (\duq)=f \mbox{ in } \Omega,$$
$$u|_{\partial \Omega}\equiv \phi.$$
\end{theorem}

\begin{remark} The real analogue of this result was proved by A.D.
Aleksandrov \cite{aleksandrov2}, and the complex one by E. Bedford
and B.A. Taylor \cite{bedford-taylor}.\end{remark}
\def\6{\partial}
\section{Generalizations to hypercomplex
manifolds.}\label{hypercomplex} Some of the definitions and
results on the quaternionic plurisubharmonic functions discussed
above were extended by M. Verbitsky and the author
\cite{alesker-verbitsky} to a more general context of so called
hypercomplex manifolds. In this section we give an overview of
these results including a geometric interpretation of quaternionic
strictly plurisubharmonic functions as (local) potentials of
HKT-metrics. The exposition follows \cite{alesker-verbitsky}.

\begin{definition}
A {\itshape hypercomplex} manifold is a smooth manifold $X$
together with a triple $(I,J,K)$ of complex structures satisfying
the usual quaternionic relations:
$$IJ=-JI=K.$$
\end{definition}
\begin{remark}
(1) We suppose here (in the opposite to much of the literature on
the subject) that the complex structures $I,J,K$ act on the
{\itshape right} on the tangent bundle $TX$ of $X$. This action
extends uniquely to the right action of the algebra $\HH$ of
quaternions on $TX$.

(2) It follows that the dimension of a hypercomplex manifold $X$
is divisible by 4.
\end{remark}

Let $(X^{4n},I,J,K)$ be a hypercomplex manifold. Let us denote by
$\Lam^{p,q}_I(X)$ the vector bundle of $(p,q)$-forms on the
complex manifold $(X,I)$. By the abuse of notation we will also
denote by the same symbol $\Lam^{p,q}_{I}(X)$ the space of
$C^\infty$-sections of this bundle.

Let
\begin{eqnarray}\label{l1}
\6\colon \Lambda_I^{p,q}(X)\to \Lambda_I^{p+1,q}(X)
\end{eqnarray}
be the usual $\6$-differential of differential forms on the
complex manifold $(X,I)$.

Set
\begin{eqnarray}\label{l2}
\6_J:=J^{-1}\circ \bar \6 \circ J.
\end{eqnarray}

\begin{claim}[\cite{verbitsky-hkt}]\label{l3}
(1)$ J\colon\Lambda_I^{p,q}(X)\to\Lam_I^{q,p}(X).$

(2) $\6_J\colon \Lambda_I^{p,q}(X)\to\Lam_I^{p+1,q}(X).$

(3) $\6\6_J=-\6_J\6$.
\end{claim}

\begin{definition}[\cite{verbitsky-hkt}]\label{l4}
Let $k=0,1,\dots,n$. A form $\ome\in \Lam^{2k,0}_I(X)$ is called
 \itshape{real} if
$$\overline{J\circ \ome}=\ome.$$
\end{definition}

We will denote the subspace of real $C^\infty$-smooth
$(2k,0)$-forms on $(X,I)$ by $\Lam^{2k,0}_{I,\RR}(X)$.
\begin{lemma}\label{l5}
Let $X$ be a hypercomplex manifold. Let $f\colon X\to \RR$ be a
smooth function. Then $\6\6_J f\in \Lam^{2,0}_{I,\RR}(X)$.
\end{lemma}

\begin{definition}
Let $\ome\in \Lam^{2,0}_{I,\RR}(X)$. Let us say that $\ome$ is
non-negative (notation: $\ome\geq 0$) if
$$\ome(Y,Y\circ J)\geq 0$$
for any (real) vector field $Y$ on the manifold $X$. Equivalently,
$\omega$ is non-negative if $\omega(Z, \bar Z \circ J)\geq 0$ for
any $(1,0)$-vector field $Z$.
\end{definition}

\begin{definition}\label{l6}
A continuous function $$h:X\to \RR$$ is called quaternionic
plurisubharmonic if $\6\6_J h$ is a non-negative (generalized)
section of $\Lambda^{2,0}_{I,\RR}(X)$.
\end{definition}
\begin{remark}
The non-negativity in the generalized sense is discussed in detail
in \cite{alesker-verbitsky}, Section 5.
\end{remark}

Let us denote by $P'(X)$ the class of continuous quaternionic
plurisubharmonic functions on $X$. Let us denote by $P''(X)$ the
subclass of functions from $P'(X)$ with the following additional
property: a function $h\in P'(X)$ belongs to $P''(X)$ if and only
if any point $x\in X$ has a neighborhood $U\ni x$ and a sequence
$\{h_N\}\subset P'(U)\cap C^2(U)$ such that $h_N\to h$ uniformly
on compact subsets of $U$. Thus $P''(X)\subset P'(X)$.

We conjecture that $P'(X)=P''(X)$. This conjecture is true when
$X$ is an open subset of $\HH^n$.

\begin{theorem}[\cite{alesker-verbitsky}, Theorem 1.10]\label{l8}
Let $X$ be a hypercomplex manifold of (real) dimension $4n$. Let
$0<k\leq n$. For any $h^{(1)},\dots,h^{(k)}\in P''(X)$ one can
define a non-negative generalized section of
$\Lambda^{2k}_{I,\RR}$ denoted by $\6\6_J h^{(1)}\wedge
\dots\wedge\6\6_J h^{(k)}$ which is uniquely characterized by the
following two properties:

(1) if $h^{(1)},\dots,h^{(k)}\in C^2(X)$ then the definition is
clear;

(2) if $\{h^{(i)}_N\}\subset P''(X)$, $h^{(i)}_N\to h^{(i)}$ as
$N\to \infty$ uniformly on compact subsets , $i=1,\dots, k$, then
$h^{(i)}\in P''(X)$ and
$$\6\6_J h^{(1)}_N\wedge \dots\wedge \6\6_J h^{(k)}_N\to \6\6_J h^{(1)}\wedge \dots\wedge \6\6_J h^{(k)}$$
in the weak topology on measures.
\end{theorem}
\begin{remark}
Theorem \ref{l8} generalizes Theorem \ref{cln2} from the flat
space $\HH^n$ to hypercomplex manifolds.
\end{remark}

Let us discuss the relations to the HKT-geometry. Let $g$ be a
Riemannian metric on a hypercomplex manifold $X$. The metric $g$
is called {\itshape quaternionic Hermitian} (or hyperhermitian) if
$g$ is invariant with respect to the group $SU(2)\subset \HH$ of
unitary quaternions.

Given a quaternionic Hermitian metric $g$ on a hypercomplex
manifold $X$, consider the differential form
$$\Ome:=\ome_J-\sqrt{-1}\ome_K$$
where $\ome_L(A,B):=g(A,B\circ L)$ for any $L\in \HH$ with
$L^2=-1$ and any real vector fields $A,B$ on $X$. It is easy to
see that $\Ome$ is a $(2,0)$-form with respect to the complex
structure $I$.

\begin{definition}\label{def-hkt-metr}
The metric $g$ on $X$ is called HKT-metric if
$$\6 \Ome =0.$$
\end{definition}
\begin{remark}
HKT-metrics on hypercomplex manifolds first were introduced by
Howe and Papadopoulos \cite{howe-papa}. Their original definition
was different but equivalent to Definition \ref{def-hkt-metr} (see
\cite{_Gra_Poon_}).
\end{remark}

Let us denote by $S_\HH(X)$ the vector bundle over $X$ such that
its fiber over a point $x\in X$ is equal to the space of
hyperhermitian forms on the tangent space $T_xX$. Consider the map
of vector bundles
$$t\colon \Lam^{2,0}_{I,\RR}(X)\to S_\HH(X)$$
defined by $t(\eta)(A,A)=\eta(A,A\circ J)$ for any (real) vector
field $A$ on $X$. Then $t$ is an isomorphism of vector bundles
(this was proved in \cite{verbitsky-hkt}).

\begin{theorem}[\cite{alesker-verbitsky}, Prop. 1.14]\label{m}
(1) Let $f$ be an infinitely smooth strictly plurisubharmonic
function on a hypercomplex manifold $(X,I,J,K)$. Then $t(\6\6_J
f)$ is an HKT-metric.

(2) Conversely assume that $g$ is an HKT-metric. Then any point
$x\in X$ has a neighborhood $U$ and an infinitely smooth strictly
plurisubharmonic function $f$ on $U$ such that $g=t(\6\6_J f)$ in
$U$.
\end{theorem}
\begin{remark}
(i) On the flat space $\HH^n$ one has for any smooth real valued
function $f$
$$t(\pt\pt_J f)=\frac{1}{4}\left(\frac{\pt^2 f}{\pt\bar q_i\pt 
q_j}\right)$$ by Proposition 4.1 of \cite{alesker-verbitsky}.

(ii) The proof of Theorem \ref{m} uses a result of Banos-Swann
\cite{banos-swann}.
\end{remark}



\end{document}